\documentclass[onecolumn]{autart}
\usepackage{amsmath}
\begin{document}
\input{amssym}
\begin{frontmatter}
\title{Symmetry classification and conservation laws for higher order Camassa-Holm equation}
\author[MN]{Mehdi Nadjafikhah}\ead{m\_nadjafikhah@iust.ac.ir},
\author[SRH]{Vahid Shirvani-Sh}\ead{v.shirvani@kiau.ac.ir}
\address[MN]{School of Mathematics, Iran University of Science and Technology, Narmak, Tehran 1684613114, Iran.}
\address[SRH]{Department of Mathematics, Islamic Azad
University, Karaj Branch, Karaj 31485-313, Iran.}
\begin{abstract}
Lie symmetry group method is applied to study for the higher order
Camassa-Holm equation. The symmetry group and its optimal system
are given. Furthermore, preliminary classification of its group
invariant solutions, symmetry reduction and nonclassical
symmeries are investigated. Finally conservation laws for the
higher order Camassa-Holm equation are presented.
\end{abstract}
\end{frontmatter}
\section{Introduction}
In the study of shallow water waves, Camassa and Holm \cite{[8]}
derived a nonlinear dispersive shallow water wave equation
\begin{eqnarray} u_t-u_{x^2t}+3uu_x=2u_xu_{x^2}+uu_{x^3},\label{eq:1.0}\end{eqnarray}
which is called Camassa-Holm equation (CH). Here $u(x,t)$ denotes
the fluid velocity at time $t$ in the $x$ direction or,
equivalently, the high of water's free surface above a flat
bottom. Eq. (1) has a bi-Hamiltonian structure \cite{[24],[10]},
and is completely integrable \cite{[8],[26]}. It has many
conservation laws \cite{[1]}. Moreover, the CH equation is a
re-expression of geodesic flow on the diffeomorphism group of the
line. Holm, Marsden and Ratiu \cite{[25]} have shown that CH
equation in $n$ dimensions describes geodesic motion on the
diffeomorphism group of $R^n$ with respect to metric given by the
$H^1$ norm or Euclidean fluid velocity. Misiolek \cite{[23]} has
shown that the CH equation represents a geodesic flow on the
Bott-Virasoro group. Kouranbaeva \cite{[22]} has shown that the
CH equation (for the case k = 0) is a geodesic spray of the weak
Riemannian metric on the diffeomorphism group of the line or the
circle obtained by the right translation of the $H^1$ inner
product over the entire group. This equation admits well known
properties and a rich literature is devoted to it.

In recent years, many researchers have been researched on the
Camassa-Holm equation. They extend the studies to the generalized
CH equation, higher order CH equations and so on. Lixin Tian,
Chunyu Shen and Danping Ding gave the optimal control of the
viscous CH equation under the boundary condition and proved the
existence and uniqueness of optimal solution to the viscous CH
equation in a short interval (see \cite{[21]}). Using geometrical
methods, higher order CH equations have been treated in
\cite{[14]}. The well-posedness of higher order CH equations were
considered in \cite{[15]}.

The formulation of the higher order Camassa-Holm equation which
was recently derived by Coclite, Holden and Karlsen in
\cite{[15]} is
\begin{eqnarray}
B_{k}(u,u):=A_{k}^{-1}C_{k}(u)-uu_{x},\nonumber\\
A_{k}(u):=\sum_{j=0}^k(-1)^j\partial_{x}^2ju,\label{eq:1.1}\\
C_k(u)=-uA_{k}(\partial_{x}u)+A_{k}(u\partial_{x}u)-2\partial_{x}uA_{k}(u).
\nonumber\end{eqnarray}
where k is a positive integer. In cases $k=0$ and $k=1$, Eq.
(\ref{eq:1.1}) becomes the inviscid Burgers equation and the
Camassa-Holm equation respectively.

In this paper we only consider the case $k=2$ of Eq.
(\ref{eq:1.1}). It also can be rewritten as
\begin{eqnarray} \Delta := u_t-u_{x^2t}+u_{x^4t}+3uu_x-2u_xu_{x^2}-uu_{x^3}+2u_{x}u_{x^4}+uu_{x^5}=0.\label{eq:1.2}\end{eqnarray}
The theory of Lie symmetry groups of differential equations was
developed by Sophus Lie \cite{[2]}, which was called classical
Lie method. Nowadays, application of Lie transformations group
theory for constructing the solutions of nonlinear partial
differential equations (PDEs) can be regarded as one of the most
active fields of research in the theory of nonlinear PDEs and
applications. Such Lie groups are invertible point transformations
of both the dependent and independent variables of the
differential equations. The symmetry group methods provide an
ultimate arsenal for analysis of differential equations and is of
great importance to understand and to construct solutions of
differential equations. Several applications of Lie groups in the
theory of differential equations were discussed in the
literature, the most important ones are: reduction of order of
ordinary differential equations, construction of invariant
solutions, mapping solutions to other solutions and the detection
of linearizing transformations. For many other applications of Lie
symmetries see \cite{[3],[4],[7]}.

The fact that symmetry reductions for many PDEs are unobtainable
by applying the classical symmetry method, motivated the creation
of several generalizations of the classical Lie group method for
symmetry reductions. The nonclassical symmetry method of
reduction was devised originally by Bluman and Cole in 1969
\cite{[11]}, to find new exact solutions of the heat equation. The
description of the method is presented in \cite{[9],[6]}. Many
authors have used the nonclassical method to solve PDEs. In
\cite{[12]} Clarkson and Mansfield have proposed an algorithm for
calculating the determining equations associated to the
nonclassical method. A new procedure for finding nonclassical
symmetries has been proposed by B\^{i}l\v{a} and Niesen in
\cite{[5]}.

Many PDEs in the applied sciences and engineering are continuity
equations which express conservation of mass, momentum, energy,
or electric charge. Such equations occur in, e.g., fluid
mechanics, particle and quantum physics, plasma physics,
elasticity, gas dynamics, electromagnetism,
magneto-hydro-dynamics, nonlinear optics, etc. In the study of
PDEs, conservation laws are important for investigating
integrability and linearization mappings and for establishing
existence and uniqueness of solutions. They are also used in the
analysis of stability and global behavior of solutions
\cite{[16],[17],[18],[19]}.

This work is organized as follows. In section 2 we recall some
results needed to construct Lie point symmetries of a given system
of differential equations. In section 3, we give the general form
of a infinitesimal generator admitted by equation (\ref{eq:1.2})
and find transformed solutions. Section 4, is devoted to the
nonclassical symmetries of the higher order CH model, symmetries
generated when a supplementary condition, the invariance surface
condition, is imposed. In Section 5, we construct the optimal
system of one-dimensional subalgebras. Lie invariants, similarity
reduced equations and differential invariants corresponding to
the infinitesimal symmetries of equation (\ref{eq:1.2}) are
obtained in section 6. Finally in last section, the conservation
laws of the equation (\ref{eq:1.2}) are obtained.
\section{Method of Lie Symmetries}
In this section, we recall the general procedure for determining
symmetries for any system of partial differential equations see
\cite{[3],[13],[4],[7]}. To begin, let us consider the general
case of a nonlinear system $E$ of partial differential equations
of order $n$ in $p$ independent and $q$ dependent variables is
given as a system of equations
\begin{eqnarray} \Delta_\nu(x,u^{(n)})=0,\;\;\;\;\; \nu=1,\cdots,l,
\label{eq:2.1} \end{eqnarray}
involving $x = (x^1,\cdots, x^p)$, $u = (u^1,\cdots,u^q)$ and the
derivatives of $u$ with respect to $x$ up to $n$, where $u^{(n)}$
represents all the derivatives of $u$ of all orders from $0$ to
$n$. We consider a one-parameter Lie group of infinitesimal
transformations acting on the independent and dependent variables
of the system (\ref{eq:2.1})
\begin{eqnarray} \tilde{x}^i &=& x^i+s \xi^i(x,u)+O(s^2), \hspace{1cm}
i=1\cdots,p,\nonumber \\[-2mm] \label{eq:2.2}\\[-2mm] \tilde{u}^j &=& u^j+s
\varphi^j(x,u)+O(s^2), \hspace{9mm} j=1\cdots,q, \nonumber
 \end{eqnarray}
where $s$ is the parameter of the transformation and $\xi^i$,
$\eta^j$ are the infinitesimals of the transformations for the
independent and dependent variables, respectively. The
infinitesimal generator ${\mathbf v}$ associated with the above
group of transformations can be written as
\begin{eqnarray}  {\mathbf v} = \sum_{i=1}^p\xi^i(x,u)\partial_{x^i} +
\sum_{\alpha=1}^q\varphi^\alpha(x,u)\partial_{u^\alpha}.
\label{eq:2.4} \end{eqnarray}
A symmetry of a differential equation is a transformation which
maps solutions of the equation to other solutions. The invariance
of the system (\ref{eq:2.1}) under the infinitesimal
transformations leads to the invariance conditions (Theorem 2.36
of \cite{[3]})
\begin{eqnarray} \textrm{Pr}^{(n)}{\mathbf
v}\big[\Delta_\nu(x,u^{(n)})\big]=0,\;\;\;\;\;
\nu=1,\cdots,l,\;\;\;\;\mbox{whenever}\;\;\;\;\;\Delta_\nu(x,u^{(n)})
=0, \label{eq:2.5} \end{eqnarray}
where $\textrm{Pr}^{(n)}$ is called the $n^{th}$ order
prolongation of the infinitesimal generator given by
\begin{eqnarray} \textrm{Pr}^{(n)}{\mathbf v}= {\mathbf
v}+\sum^q_{\alpha=1}\sum_J
\varphi^J_\alpha(x,u^{(n)})\partial_{u^\alpha_J},\label{eq:2.6}
\end{eqnarray}
where $J=(j_1,\cdots,j_k)$, $1\leq j_k\leq p$, $1\leq k\leq n$ and
the sum is over all $J$'s of order $0<\# J\leq n$. If $\#J=k$, the
coefficient $\varphi_J^\alpha$ of $\partial_{u_J^\alpha}$ will
only depend on $k$-th and lower order derivatives of $u$, and
\begin{eqnarray} \varphi_\alpha^J(x,u^{(n)})=D_J(\varphi_\alpha - \sum_{i=1}^p
\xi^iu_i^\alpha) + \sum_{i=1}^p\xi^iu^\alpha_{J,i}, \label{eq:2.7}
\end{eqnarray}
where $u_i^\alpha:=\partial u^\alpha/\partial x^i$ and
$u_{J,i}^\alpha := \partial u_J^\alpha/\partial x^i$.

\medskip One of the most important properties of these
infinitesimal symmetries is that they form a Lie algebra under the
usual Lie bracket.
\section{Lie symmetries for the higher order CH equation  }
We consider the one parameter Lie group of infinitesimal
transformations on $(x^1=x,x^2=t,u^1=u)$,
\begin{eqnarray} \tilde{x} &=& x+s\xi(x,t,u)+O(s^2),\nonumber\\
\tilde{t} &=& x+s\eta(x,t,u)+O(s^2),\label{eq:3.1}\\
\tilde{u} &=& x+s\varphi(x,t,u)+O(s^2),\nonumber \end{eqnarray}
where $s$ is the group parameter and $\xi^1=\xi$, $\xi^2=\eta$ and
$\varphi^1=\varphi$ are the infinitesimals of the transformations
for the independent and dependent variables, respectively. The
associated vector field is of the form:
\begin{eqnarray} {\mathbf
v}=\xi(x,t,u)\partial_x+\eta(x,t,u)\partial_t+\varphi(x,t,u)\partial_u.
\label{eq:3.2}\end{eqnarray}
and, by (\ref{eq:2.6}) its fifth prolongation is
\begin{eqnarray} \textrm{Pr}^{(5)}{\mathbf v} &=& {\mathbf v}+
\varphi^x\,\partial_{u_x}+\varphi^t\,\partial_{u_t}+\varphi^{x^2}\,\partial_{u_{x^2}}
+\varphi^{xt}\,\partial_{u_{xt}}+\varphi^{t^2}\,\partial_{u_{t^2}} \nonumber\\
&&
+\varphi^{x^3}\,\partial_{u_{x^3}}+\varphi^{x^2t}\,\partial_{u_{x^2t}}+\varphi^{xt^2}\,\partial_{u_{xt^2}}+\varphi^{t^3}\,\partial_{u_{t^3}}
+\cdots+\varphi^{xt^4}\,\partial_{u_{xt^4}}+\varphi^{t^5}\,\partial_{u_{t^5}}.
\label{eq:3.2-1}\end{eqnarray}
where, for instance by (\ref{eq:2.7}) we have
\begin{eqnarray}
\varphi^x&=&D_x(\varphi-\xi\,u_x-\eta\,u_t)+\xi\,u_{x^2}+\eta\,u_{xt},\nonumber\\
\varphi^t&=&D_t(\varphi-\xi\,u_x-\eta\,u_{t})+\xi\,u_{xt}+\eta\,u_{t^2},\nonumber\\
&& \vdots \label{eq:3.3} \\
\varphi^{t^5}&=&D^5_{x}(\varphi-\xi\,u_x-\eta\,u_t)+\xi\,u_{x^5t}+\eta\,u_{t^5},\nonumber
\end{eqnarray}
where $D_x$ and $D_t$ are the total derivatives with respect to
$x$ and $t$ respectively.
By (\ref{eq:2.5}) the vector field ${\mathbf v}$ generates a one
parameter symmetry group of the Eq. (\ref{eq:1.2}) if and only if
\begin{eqnarray}  \textrm{Pr}^{(5)}{\mathbf
v}[u_t-u_{x^2t}+u_{x^4t}+3uu_x-2u_xu_{x^2}-uu_{x^3}+2u_{x}u_{x^4}+uu_{x^5}]=0,\,\,\,
\mbox{whenever} \hspace{.5cm} \Delta=0.
\label{eq:3.3-1}\end{eqnarray}
The condition (\ref{eq:3.3-1}) is equivalent to
\begin{eqnarray}
(3u_{x}-u_{x^3}+u_{x^5})\varphi+\varphi^t+(3u-2u_{x^2}+2u_{x^4})\varphi^x-2u_{x}\varphi^{x^2}-
u\varphi^{x^3}-\varphi^{x^2t}+2u_{x}\varphi^{x^4}+u\varphi^{x^5}+\varphi^{x^4t}=0,\nonumber\\
\label{eq:3.4}\\
 \mbox{whenever}\hspace{.3cm}u_t-u_{x^2t}+u_{x^4t}+3uu_x-2u_xu_{x^2}-uu_{x^3}+2u_{x}u_{x^4}+uu_{x^5}=0.\nonumber
\end{eqnarray}
Substituting (\ref{eq:3.3}) into (\ref{eq:3.4}), and equating the
coefficients of the various monomials in partial derivatives with
respect to $x$ and various power of $u$, we can find the
determining equations for the symmetry group of the Eq.
(\ref{eq:1.2}). Solving this equations, we get the following
forms of the coefficient functions
\begin{eqnarray}
 \xi=c_3, \quad\quad\quad\quad \eta=c_1t+c_2, \quad\quad\quad\quad  \varphi=-c_1u. \label{eq:3.5}\end{eqnarray}
where $c_1$, $c_2$ and $c_3$ are arbitrary constant. Thus, the Lie
algebra $\goth g$ of infinitesimal symmetry of the Eq.
(\ref{eq:1.2}) is spanned bye the three vector fields
\begin{eqnarray} \textbf{v}_1=\partial_x,\quad\quad\quad\quad \textbf{v}_2=\partial_t,\quad\quad\quad\quad
\textbf{v}_3=t\partial_t-u\partial_u. \label{eq:3.6}
\end{eqnarray}
The commutation relations between these vector fields are given
in the Table 1. The Lie algebra $\goth g$ is solvable, because if
${\goth g}^{(1)}=\langle
\textbf{v}_i,[\textbf{v}_i,\textbf{v}_j]\rangle=[\goth g,\goth
g]$, we have ${\goth g}^{(1)}=\langle \textbf{v}_1,\textbf{v}_2,
\textbf{v}_3\rangle$, and ${\goth g}^{(2)}=[{\goth g}^{(1)},{\goth
g}^{(1)}]=\langle\textbf{v}_2\rangle$, so, we have a chain of
ideals ${\goth g}^{(1)}\supset {\goth g}^{(2)}\supset {0}$.
\begin{table}[h] \label{Tab:1}
\small \caption{\textit{The commutator table}}
\begin{tabular}{cccc} \hline
  $[{\mathbf v}_{i},{\mathbf v}_{j}]$ & ${\mathbf v}_1$ & ${\mathbf v}_2$ & ${\mathbf v}_3$  \\ \hline
  ${\mathbf v}_1$ & 0 & 0 & 0  \\
  ${\mathbf v}_2$ & 0 & 0 & ${\mathbf v}_2$ \\
  ${\mathbf v}_3$ & 0 & $-{\mathbf v}_2$ & 0 \\\hline
 \end{tabular}
\end{table}

\medskip To obtain the group transformation which is generated by the
infinitesimal generators $\textbf{v}_i$ for $i=1,2,3$ we need to
solve the three systems of first order ordinary differential
equations
\begin{eqnarray*}
 \frac{d\tilde{x}(s)}{ds} &=&
\xi_i(\tilde{x}(s),\tilde{t}(s),\tilde{u}(s)), \quad
\tilde{x}(0)=x, \nonumber\\
 \frac{d\tilde{t}(s)}{ds} &=&
\eta_i(\tilde{x}(s),\tilde{t}(s),\tilde{u}(s)), \quad
\tilde{t}(0)=t, \qquad i=1,2,3 \label{eq:3.7}\\
 \frac{d\tilde{u}(s)}{ds} &=&
\varphi_i(\tilde{x}(s),\tilde{t}(s),\tilde{u}(s)), \quad
\tilde{u}(0)=u. \nonumber \end{eqnarray*}
Exponentiating the infinitesimal symmetries of equation
(\ref{eq:1.2}), we get the one-parameter groups $G_i(s)$ generated
by $\textbf{v}_i$ for $i=1,2,3$
\begin{eqnarray}
G_1:(t,x,u) & \longmapsto & (x+s,t,u),\nonumber\\
G_2:(t,x,u) & \longmapsto & (x,t+s,u),\label{eq:3.8} \\
G_3:(t,x,u) &\longmapsto&(x,{\rm {e}}^{s}t,{\rm
{e}}^{-s}u).\nonumber \end{eqnarray}
Consequently,

\textbf{Theorem 3.1} {\it If $u=f(x,t)$ is a solution of higher
order CH equation, so are the functions
\begin{eqnarray} G_1(s)\cdot f(x,t)&=&f(x-s,t),\nonumber\\
G_2(s)\cdot f(x,t)&=&f(x,t-s),\label{eq:3.9} \\
G_3(s)\cdot f(x,t)&=&f(x,t{\rm {e}}^{-s}){\rm {e}}^{-s}.\nonumber
\end{eqnarray} }
\section{Nonclassical symmetries for the higher order CH equation }
In this section we would like to apply the nonclassical method to
the higher order CH equation. The graph of a solution
\begin{eqnarray}
u^{\alpha}=f^{\alpha}(x_{1},\cdots,x_{p}),\hspace{1.5cm}\alpha=1,\cdots,q\label{eq:4.1}\end{eqnarray}
to the system (\ref{eq:2.1}) defines an p-dimensional submanifold
$\Gamma_{f}\subset {\mathbf{R}}^{p}\times {\mathbf{R}}^{q}$ of the
space of independent and dependent variables. The solution will be
invariant under the one-parameter subgroup generated by vector
(\ref{eq:2.4}) if and only if $\Gamma_{f}$ is an invariant
submanifold of this group. By applying the well known criterion
of invariance of a submanifold under a vector field we get that
(\ref{eq:4.1}) is invariant under vector (\ref{eq:2.4}) if and
only if $f$ satisfies the first order system $E_{Q}$ of partial
differential equations
\begin{eqnarray}
Q^{\alpha}(x,u,u^{(1)})=\varphi^{(\alpha)}(x,u)-\sum_{i=1}^p\xi^i(x,u)u_i^\alpha=0,\hspace{1cm}\alpha=1,\cdots,q\label{eq:4.2}\end{eqnarray}
known as the invariant surface conditions. The q-tuple
$Q=(Q^1,\cdots,Q^q)$ is known as the characteristic of the vector
field (\ref{eq:2.4}). In what follows, the n-th prolongation of
the invariant surface conditions (\ref{eq:4.2}) will be denoted
by $E_Q^{(n)}$, which is a n-th order system of partial
differential equations obtained by appending to (\ref{eq:4.2}) its
partial derivatives with respect to the independent variables of
orders $j\leq n-1$.

For the system (\ref{eq:2.1}), (\ref{eq:4.2}) to be compatible,
the n-th prolongation $\textrm{Pr}^{(n)}{\mathbf v}$ of the
vector field ${\mathbf v}$ must be tangent to the intersection
$E\cap E_{Q}^{(n)}$
\begin{eqnarray} \textrm{Pr}^{(n)}{\mathbf v}(\Delta_\nu)|_{E\cap
E_{Q}^{(n)}}=0, \hspace
{1cm}\nu=1,\cdots,l.\label{eq:4.3}\end{eqnarray}
If the equations (\ref{eq:4.3}) are satisfied, then the vector
field (\ref{eq:2.4}) is called a nonclassical infinitesimal
symmetry of the system (\ref{eq:2.1}). The relations
(\ref{eq:4.3}) are generalizations of the relations
(\ref{eq:2.5}) for the vector fields of the infinitesimal
classical symmetries. A similar procedure is applicable to the
case of the nonclassical infinitesimal symmetries with an evident
difference that in general one has fewer determining equations
than in the classical case. Therefore, we expect that nonclassical
symmetries are much more numerous than classical ones, since any
classical symmetry is clearly a nonclassical one. The important
feature of determining equations for nonclassical symmetries is
that they are nonlinear, this implies that the space of
nonclassical symmetries does not, in general, form a vector
space. For more theoretical background see \cite{[20],[5]}.

If we assume that the coefficient of $\partial_t$ of the vector
field (\ref{eq:2.4}) does not identically equal zero, then for
the vector field
\begin{eqnarray} {\mathbf
v}=\xi(x,t,u)\partial_x+\partial_t+\varphi(x,t,u)\partial_u
\label{eq:4.5}\end{eqnarray}
the invariant surface conditions are
\begin{eqnarray} u_t+\xi u_x=\varphi,
\label{eq:4.6}\end{eqnarray}
Calculating equations (\ref{eq:4.3}) and inserting $\varphi$ from
(\ref{eq:4.6}) in to it, we can find the determining equations by
equating the coefficients of the various monomials in partial
derivatives with respect to $x$ and various power of $u$. Solving
this equations, we get $\xi=c$ and $\varphi=0$, where $c$ is
arbitrary constant.

Now assume that the coefficient of $\partial_t$ in (\ref{eq:4.5})
equals zero and try to find the infinitesimal nonclassical
symmetries of the form
\begin{eqnarray} {\mathbf
v}=\partial_x+\varphi(x,t,u)\partial_u
\label{eq:4.9}\end{eqnarray}
for which the invariant surface conditions is $u_x=\varphi$.
Similar the previous case, we can find determining equations.
Solving this equations, we get $\varphi=0$. This means that no
supplementary symmetries, of non-classical type, are specific for
our models.
\section{Optimal system for the higher order CH equation }
In general, to each s-parameter subgroup $H$ of the full symmetry
group $G$ of a system of differential equations in $p>s$
independent variables, there will correspond a family of
group-invariant solutions. Since there are almost always an
infinite number of such subgroups, it is not usually feasible to
list all possible group-invariant solutions to the system. We
need an effective, systematic means of classifying these
solutions, leading to an "optimal system" of group-invariant
solutions from which every other such solution can be derived.

\textbf{Definition 5.1} Let $G$ be a Lie group with Lie algebra
$\goth g$. An optimal system of $s-$parameter subgroups is a list
of conjugacy inequivalent $s-$parameter subalgebras with the
property that any other subgroup is conjugate to precisely one
subgroup in the list. Similarly, a list of $s-$parameter
subalgebras forms an optimal system if every $s-$parameter
subalgebra of $\goth g$ is equivalent to a unique member of the
list under some element of the adjoint representation:
$\overline{\goth h}={\mathrm Ad}(g({\goth h}))$, $g\in
G$.\cite{[3]}

\textbf{Theorem 5.2} {\it Let $H$ and $\overline{H}$ be connected
s-dimensional Lie subgroups of the Lie group $G$ with
corresponding Lie subalgebras $\goth h$ and $\overline{\goth h}$
of the Lie algebra $\goth g$ of $G$. Then
$\overline{H}$=$gHg^{-1}$ are conjugate subgroups if and only if
$\overline{\goth h}={\mathrm Ad}(g({\goth h}))$ are conjugate
subalgebras.}(Proposition 3.7 of \cite{[3]})

\medskip By theorem (5.2), the problem of finding an optimal system of
subgroups is equivalent to that of finding an optimal system of
subalgebras. For one-dimensional subalgebras, this classification
problem is essentially the same as the problem of classifying the
orbits of the adjoint representation, since each one-dimensional
subalgebra is determined by nonzero vector in $\goth g$. This
problem is attacked by the na\"{\i}ve approach of taking a general
element ${\mathbf V}$ in $\goth g$ and subjecting it to various
adjoint transformation so as to "simplify" it as much as
possible. Thus we will deal with the construction of the optimal
system of subalgebras of $\goth g$. To compute the adjoint
representation, we use the Lie series
\begin{eqnarray} {\mathrm Ad}(\exp(\varepsilon{\mathbf v}_i){\mathbf v}_j) =
{\mathbf v}_j-\varepsilon[{\mathbf v}_i,{\mathbf
v}_j]+\frac{\varepsilon^2}{2}[{\mathbf v}_i,[{\mathbf
v}_i,{\mathbf v}_j]]-\cdots,\label{eq:5.1} \end{eqnarray}
where $[{\mathbf v}_i,{\mathbf v}_j]$ is the commutator for the
Lie algebra, $\varepsilon$ is a parameter, and $i,j=1,2,3$. Then
we have the Table 2.

\begin{table}[h] \label{Tab:2}
\small \caption{\textit{Adjoint representation table }}
\begin{tabular}{cccc} \hline
  $Ad(\exp(\varepsilon{\mathbf v}_i){\mathbf v}_j)$ & ${\mathbf v}_1$ & ${\mathbf v}_2$ & ${\mathbf v}_3$ \\\hline
  ${\mathbf v}_1$ & ${\mathbf v}_1$ & ${\mathbf v}_2$ & ${\mathbf v}_3$  \\
  ${\mathbf v}_2$ & ${\mathbf v}_1$ & ${\mathbf v}_2$ & ${\mathbf v}_3-\varepsilon{\mathbf v}_2$ \\
  ${\mathbf v}_3$ & ${\mathbf v}_1$ & ${\mathbf v}_2+\varepsilon{\mathbf v}_3$ & ${\mathbf v}_3$  \\ \hline
\end{tabular}
\end{table}

\textbf{Theorem 5.3} {\it An optimal system of one-dimensional Lie
algebras of the higher order CH equation is provided by\\
\hspace{2cm}(1) $\;\alpha\textbf{v}_1+\textbf{v}_3$,
\quad\quad\quad (2) $\;\beta\textbf{v}_1+\textbf{v}_2$ }

\textbf{proof}: Consider the symmetry algebra $\goth g$ of the Eq.
(\ref{eq:1.2}) whose adjoint representation was determined in
table 2 and let $F^s_i:{\goth g}\to{\goth g}$ defined by
${\mathbf v}\mapsto\mathrm{Ad}(\exp(\varepsilon{\mathbf
v}_i){\mathbf v})$ is a linear map, for $i=1,2,3$. The matrices
$M^\varepsilon_i$ of $F^\varepsilon_i$, $i=1,2,3$, with respect to
basis $\{{\mathbf v}_1,{\mathbf v}_2,{\mathbf v}_3\}$ are
\begin{eqnarray*}
M^\varepsilon_1= \left( \begin {array}{ccc}
1&0&0\\
0&1&0\\
0&0&1\end {array} \right),\,\,\quad\quad
 M^\varepsilon_2= \left( \begin {array}{ccc}
1&0&0\\
0&1&\varepsilon\\
0&0&1\end {array} \right),\,\,\quad\quad
 M^\varepsilon_3= \left( \begin {array}{ccc}
1&0&0\\
0&{{\rm e}^{{-\varepsilon}}}&0\\
0&0&1
\end {array} \right).\,\,
\end{eqnarray*}
Let ${\mathbf V}=\sum_{i=1}^3a_i{\mathbf v}_i$ is a nonzero
vector field in $\goth g$. We will simplify as many of the
coefficients $a_{i}$ as possible by acting these matrices on a
vector field ${\mathbf V}$ alternatively.

Suppose first that $a_{3}\neq0$, scaling ${\mathbf V}$ if
necessary we can assume that $a_{3}=1$, then we can make the
coefficients of ${\mathbf v}_2$ vanish by $M^\varepsilon_2$, and
${\mathbf V}$ reduced to case 1.

If $a_{3}=0$ and $a_{2}\neq0$, then we can not make vanish the
coefficients of ${\mathbf v}_1$ and ${\mathbf v}_2$ by acting any
matrices $M^\varepsilon_i$. Scaling ${\mathbf V}$ if necessary,
we can assume that $a_2=1$ and ${\mathbf V}$ reduced to case
2.\hspace{7cm} $\Box$
\section{Symmetry reduction and  differential invariants for the higher order CH equation }
Lie-group method is applicable to both linear and non-linear
partial differential equations, which leads to similarity
variables that may be used to reduce the number of independent
variables in partial differential equations. By determining the
transformation group under which a given partial differential
equation is invariant, we can obtain information about the
invariants and symmetries of that equation.

Symmetry group method will be applied to the (\ref{eq:1.2}) to be
connected directly to some order differential equations. To do
this, a particular linear combinations of infinitesimals are
considered and their corresponding invariants are determined. The
equation (\ref{eq:1.2}) is expressed in the coordinates
$(x,t,u)$, so to reduce this equation is to search for its form
in specific coordinates. Those coordinates will be constructed by
searching for independent invariants $(y,v)$ corresponding to the
infinitesimal generator. So using the chain rule, the expression
of the equation in the new coordinate allows us to the reduced
equation. Here we will obtain some invariant solutions with
respect to symmetries. First we obtain the similarity variables
for each term of the Lie algebra $\goth g$, then we use this
method to reduced the PDE and find the invariant solutions.

We can now compute the invariants associated with the symmetry
operators, they can be obtained by integrating the characteristic
equations. For example for the operator
$\alpha\textbf{v}_1+\textbf{v}_2=\alpha\partial_x+\partial_t$
characteristic equation is
\begin{eqnarray} \frac{dx}{\alpha}=\frac{dt}{1}=\frac{du}{0}. \label{eq:6.1} \end{eqnarray}
The corresponding invariants are $y=x-\alpha\,t$, $v=u$
therefore, a solution of our equation in this case is $u=v(y)$.
The derivatives of $u$ are given in terms of $y$ and $v$ as
\begin{eqnarray} u_t=-\alpha\,v_{y},\quad u_{x^2t}=-\alpha\,v_{y^3},\quad
u_{x^4t}=-\alpha\,v_{y^5},\quad u_x=v_{y},\quad
u_{x^2}=v_{y^2},\cdots,
u_{x^5}=v_{y^5}.\label{eq:6.2}\end{eqnarray}
Substituting (\ref{eq:6.2}) into the Eq. (\ref{eq:1.2}), we obtain
the ordinary differential equation
\begin{eqnarray} -\alpha\,v_{y}+\alpha\,v_{y^3}-\alpha\,v_{y^5}
+3\,v\,v _{y}-2\,v_{y} v_{y^2}-v\,v
_{y^3}+2\,v_{y}v_{y^4}+v\,v_{y^5}=0 .\label{eq:6.3}\end{eqnarray}
All results are coming in the tables 3 and 4.
\begin{table}[h] \label{Tab:3}
\small \caption{\textit{Reduction of Eq. (\ref{eq:1.2}) }} \small
\begin{tabular}{cccc} \hline
  operator & $y$ & $v$ & $u$ \\\hline
  ${\mathbf v}_1$ & $t$ & $u$ & $v(y)$ \\
  ${\mathbf v}_2$ & $x$ & $u$ & $v(y)$\\
  ${\mathbf v}_3$ & $x$ & $t\,u$ & $\frac{1}{t}v(y)$ \\
  $\alpha\,{\mathbf v}_1+{\mathbf v}_3$ & $x-\log(t)$ & $t\,u$ & $\frac{1}{t}v(y)$  \\
  $\alpha\,{\mathbf v}_1+{\mathbf v}_2$ & $x-\alpha\,t$ & $u$ & $v(y)$ \\
  \hline
\end{tabular}
\end{table}

\begin{table}[h] \label{Tab:4}
\small \caption{\textit{Reduced equations corresponding to
infinitesimal symmetries }}

\small \begin{tabular}{cc} \hline
  operator & similarity reduced equations \\\hline
  ${\mathbf v}_1$ & $v_{y}=0$ \\
  ${\mathbf v}_2$ & $3vv_{y}-2v_{y}v_{y^2}-vv_{y^3}+2v_{y}v_{y^4}+vv_{y^5}=0$\\
  ${\mathbf v}_3$ & $-v+v_{y^2}-v_{y^4}+3vv_{y}-2v_{y}v_{y^2}-vv_{y^3}+2v_{y}v_{y^4}+vv_{y^5}=0$ \\
  $\alpha\,{\mathbf v}_1+{\mathbf v}_2$ & $-\alpha\,v_{y}+\alpha\,v_{y^3}-\alpha\,v_{y^5}+3\,v\,v _{y}-2\,v_{y} v_{y^2}-v\,v_{y^3}+2\,v_{y}v_{y^4}+v\,v_{y^5}=0$ \\
  $\alpha\,{\mathbf v}_1+{\mathbf v}_3$ & $-v-\alpha\,v_{y}+v_{y^2}+\alpha\,v_{y^3}-v_{y^4}-\alpha\,v_{y^5}+3vv_{y}-2v_{y}v_{y^2}-vv_{y^3}+2v_{y}v_{y^4}+vv_{y^5}=0$\\
  \hline
\end{tabular}
\end{table}

Differential invariants help us to find general systems of
differential equations which admit a prescribed symmetry group.
One say, if G is a symmetry group for a system of PDEs with
functionally differential invariants, then, the system can be
rewritten in terms of differential invariants. For finding the
differential invariants of the equation (\ref{eq:1.2}) up to
order 2, we should solve the following systems of PDEs:
\begin{eqnarray} {\frac {\partial I}{\partial x}},\quad \quad{\frac
{\partial I}{\partial t}}, \quad \quad t\,{\frac {\partial
I}{\partial t}}-u\,{\frac {\partial I}{\partial u}},
\label{eq:6.4}\end{eqnarray}

where $I$ is a smooth function of $(x,t,u)$,
\begin{eqnarray}{\frac {\partial I_1}{\partial x}},\quad \quad {\frac {\partial
I_1}{\partial t}},\quad \quad t\,{\frac {\partial I_1}{\partial
t}}-u\,{\frac {\partial I_1}{\partial u}}-u_{x}{\frac {\partial
I_1}{\partial u_{x}}}-2u_{t}\,{\frac {\partial I_1}{\partial
u_{t}}}, \label{eq:6.5}\end{eqnarray}
where $I_{1}$ is a smooth function of $(x,t,u,u_{x},u_{t})$,
\begin{eqnarray} {\frac {\partial I_2}{\partial x}},\quad \quad {\frac {\partial
I_2}{\partial t}},\quad \quad t\,{\frac {\partial I_2}{\partial
t}}-u\,{\frac {\partial I_2}{\partial u}}-u_{x}{\frac {\partial
I_1}{\partial u_{x}}}-2u_{t}\,{\frac {\partial I_2}{\partial
u_{t}}}-u_{x^2}\,{\frac {\partial I_2}{\partial
u_{x^2}}}-2\,u_{xt}\,{\frac {\partial I_2}{\partial
u_{xt}}}-3\,u_{t^2}\,{\frac {\partial I_2}{\partial u_{t^2}}},
\label{eq:6.6}\end{eqnarray}
where $I_{2}$ is a smooth function of
$(x,t,u,u_{x},u_{t},u_{xx},u_{xt},u_{tt})$. The solutions of PDEs
systems (\ref{eq:6.4}),(\ref{eq:6.5}) and (\ref{eq:6.6}) coming
in table 5, where * and ** are refer to ordinary invariants and
first order differential invariants respectively.
\begin{table}[h] \label{Tab:5}
\small \caption{\textit{differential invariants }}

\small \begin{tabular}{cccc} \hline
  vector field & ordinary invariant  & 1st order & 2nd order \\\hline
  ${\mathbf v}_1$ & $t,u$ & $*,u_{x},u_{t}$ & $*,**,u_{xx},u_{xt},u_{tt}$ \\
  ${\mathbf v}_2$ & $x,u$ & $*,u_{x},u_{t}$ & $*,**,u_{xx},u_{xt},u_{tt}$\\
  ${\mathbf v}_3$ & $x,t\,u$ & $*,t\,u_{x},t^2\,u_{t}$ & $*,**,t\,u_{xx},t^2\,u_{xt},t^3\,u_{tt}$ \\
  \hline
\end{tabular}
\end{table}

\section{Conservation laws for the higher order CH equation}
Many methods for dealing with the conservation laws are derived,
such as the method based on the Noether's theorem, the multiplier
method, by the relationship between the conserved vector of a PDE
and the Lie-Bäcklund symmetry generators of the PDE, the direct
method, etc.\cite{[3],[16],[17],[18]}.

Now, we derive the conservation laws from the multiplier method.

\textbf{Definition 8.1} A local conservation law of the PDE system
(\ref{eq:2.1}) is a divergence expression
\begin{eqnarray} D_i\Phi^i[u]=D_1\Phi^1[u]+\cdots+D_n\Phi^n[u]=0
\label{eq:8.1}\end{eqnarray}
holding for all solutions of the system (\ref{eq:2.1}). In
(\ref{eq:8.1}),
$\Phi^i[u]=\Phi^i(x,u,\partial_u,\cdots,\partial^r_u)$,
$i=1,\cdots,n$, are called fluxes of the conservation law, and
the highest-order derivative $(r)$ present in the fluxes
$\Phi^i[u]$ is called the order of a conservation law. \cite{[17]}

\medskip In particular, a set of multipliers
$\{\Lambda_\nu[U]\}^l_{\nu=1}=\{\Lambda_\nu(x,U,\partial_U,\cdots,\partial^r_U)\}^l_{\nu=1}$
yields a divergence expression for the system
$\Delta_\nu(x,u^{(n)})$ such that if the identity
\begin{eqnarray} \Lambda_\nu[U]\Delta_\nu[U] \equiv D_i\Phi^i[U]
\label{eq:8.2}\end{eqnarray}
holds identically for arbitrary functions $U(x)$. Then on the
solutions $U(x)=u(x)$ of the system (\ref{eq:2.1}), if
$\Lambda_\nu[U]$ is non-singular, one has local conservation law
$\Lambda_\nu[u]\Delta_\nu[u]=D_i\Phi^i[u]=0$.

\textbf{Definition 8.2} The Euler operator with respect to $U^j$
is the operator defined by
\begin{eqnarray}
E_{U^j}=\frac{\partial}{\partial{U^j}}-D_i\frac{\partial}{\partial{U^j}}+\cdots+(-1)^sD_{i_1}\cdots
D_{i_s}\frac{\partial}{\partial{U^j_{i_1\cdots i_s}}}+\cdots
\label{eq:8.3}\end{eqnarray}
for $j=1,\cdots,q$.\cite{[17]}

\textbf{Theorem 8.3} {\it The equations
$E_{U^j}F(x,U,\partial_U,\cdots,\partial^s_U)\equiv0$,
$j=1,\cdots,q$ hold for arbitrary $U(x)$ if and only if
$F(x,U,\partial_U,\cdots,\partial^s_U)\equiv
D_i\Psi^i(x,U,\partial_U,\cdots,\partial^{s-1}_U)$ holds for some
functions $\Psi^i(x,U,\partial_U,\cdots,\partial^{s-1}_U)$,\,
$i=1,\cdots q$.} (Theorem 1.3.2, \cite{[17]})

\textbf{Theorem 8.4} {\it A set of non-singular local multipliers
$\{\Lambda_\nu(x,U,\partial_U,\cdots,\partial^r_U)\}^l_{\nu=1}$
yields a local conservation law for the system
$\Delta_\nu(x,u^{(n)})$ if and only if the set of identities
\begin{eqnarray}
E_{U^j}(\Lambda_\nu(x,U,\partial_U,\cdots,\partial^r_U)\Delta_\nu(x,u^{(n)}))\equiv
0, \,j=1,\cdots q,\label{eq:8.4}\end{eqnarray} holds for
arbitrary functions $U(x)$.} (Theorem 1.3.3, \cite{[17]})

\medskip The set of equations (\ref{eq:8.4}) yields the set of linear
determining equations to find all sets of local conservation law
multipliers of the system (\ref{eq:2.1}). Now, we seek all local
conservation law multipliers of the form $\Lambda=\xi(x,t,u)$ of
the equation (\ref{eq:1.2}). The determining equations
(\ref{eq:8.4}) become
\begin{eqnarray}
E_U[\xi(x,t,U)(U_t-U_{x^2t}+U_{x^4t}+3UU_x-2U_xU_{x^2}-UU_{x^3}+2U_{x}U_{x^4}+UU_{x^5})]\equiv
0,\label{eq:8.5}\end{eqnarray}
where $U(x,t)$ are arbitrary function. Equation (\ref{eq:8.5})
split with respect to third order derivatives of $U$ to yield the
determining PDE system whose solutions are the sets of local
multipliers of all nontrivial local conservation laws of the
higher order CH equation.

The solution of the determining system (\ref{eq:8.5}) given by
\begin{eqnarray} \xi=c_1\,U+c_2,\label{eq:8.6}\end{eqnarray}
where $c_1$ and $c_2$ are arbitrary constants. So local
multipliers given by
\begin{eqnarray} 1)\; \xi=1, \hspace{1cm} 2)\; \xi=U,\label{eq:8.7}\end{eqnarray}
Each of the local multipliers $\xi$ determines a nontrivial local
conservation law $D_t\Psi+D_x\Phi=0$ with the characteristic form
\begin{eqnarray} D_t\Psi+D_x\Phi \equiv
\xi(U_t-U_{x^2t}+U_{x^4t}+3UU_x-2U_xU_{x^2}-UU_{x^3}+2U_{x}U_{x^4}+UU_{x^5}),\label{eq:8.8}\end{eqnarray}
To calculate the conserved quantities $\Psi$ and $\Phi$, we need
to invert the total divergence operator. This requires the
integration (by parts) of an expression in multi-dimensions
involving arbitrary functions and its derivatives, which is a
difficult and cumbersome task. The homotopy operator \cite{[19]}
is a powerful algorithmic tool (explicit formula) that originates
from homological algebra and variational bi-complexes.

\textbf{Definition 8.5} The 2-dimensional homotopy operator is a
vector operator with two components,
$\Big({\textit{H}}^{(x)}_{u(x,t)}\textit{f},{\textit{H}}^{(t)}_{u(x,t)}\textit{f}\Big)$,
where
\begin{eqnarray}
{\textit{H}}^{(x)}_{u(x,t)}\textit{f}=\int^1_0\Big(\sum^q_{j=1}\textit{I}^{(x)}_{u^j}\textit{f}\Big)[\lambda
u]\frac{d\lambda}{\lambda}\;\;\;\hspace{.2cm} and\;\;\;
\hspace{.2cm}
{\textit{H}}^{(t)}_{u(x,t)}\textit{f}=\int^1_0\Big(\sum^q_{j=1}\textit{I}^{(t)}_{u^j}\textit{f}\Big)[\lambda
u]\frac{d\lambda}{\lambda}.\label{eq:8.9}\end{eqnarray}
The x-integrand, $\textit{I}^{(x)}_{u^j_{(x,t)}}\textit{f}$, is
given by
\begin{eqnarray}
\textit{I}^{(x)}_{u^j}\textit{f}=\sum^{M^j_1}_{k_1=1}\sum^{M^j_2}_{k_2=0}\Big(\sum^{k_1-1}_{i_1=0}\sum^{k_2}_{i_2=0}
B^{(x)}u^j_{x^{i_1}t^{i_2}}(-D_x)^{k_1-i_1-1}(-D_t)^{k_2-i_2}\Big)\frac{\partial\textit{f}}{\partial
u^j_{x^{k_1}t^{k_2}}},\label{eq:8.10}\end{eqnarray}
where $M^j_1$, $M^j_2$ are the order of $\textit{f}$ in $u$ to
$x$ and $t$ respectively and combinatorial coefficient

\begin{eqnarray*}
B^{(x)}=B(i_1,i_2,k_1,k_2)=\frac{\left(\!\!\!\begin{array}{c}i_1+i_2\\i_1
\end{array}\!\!\!\right)\left(\!\!\!\begin{array}{c}k_1+k_2-i_1-i_2-1\\k_1-i_1-1
\end{array}\!\!\!\right)}{\left(\!\!\!\begin{array}{c}k_1+k_2\\k_1
\end{array}\!\!\!\right)}.
\label{eq:8.11}\end{eqnarray*}
Similarly, the t-integrand,
$\textit{I}^{(t)}_{u^j_{(x,t)}}\textit{f}$, is defined as
\begin{eqnarray}
\textit{I}^{(t)}_{u^j}\textit{f}=\sum^{M^j_1}_{k_1=0}\sum^{M^j_2}_{k_2=1}\Big(\sum^{k_1}_{i_1=0}\sum^{k_2-1}_{i_2=0}
B^{(t)}u^j_{x^{i_1}t^{i_2}}(-D_x)^{k_1-i_1}(-D_t)^{k_2-i_2-1}\Big)\frac{\partial\textit{f}}{\partial
u^j_{x^{k_1}t^{k_2}}},\label{eq:8.12}\end{eqnarray}
where $B^{(t)}(i_2,i_1,k_2,k_1)$.

\medskip We apply homotopy operator to find conserved
quantities $\Psi$ and $\Phi$ which yield of multiplier $\xi=1$.
We have
\begin{eqnarray}
\textit{f}=U_t-U_{x^2t}+U_{x^4t}+3UU_x-2U_xU_{x^2}-UU_{x^3}+2U_{x}U_{x^4}+UU_{x^5},\label{eq:8.13}\end{eqnarray}
the integrands (\ref{eq:8.10}) and (\ref{eq:8.12}) are
\begin{eqnarray}
\begin{array}{lcl}
\textit{I}^{(x)}_{u^j}\textit{f}=3\,u^{2}-u_{x}^2-\frac{2}{3}\,u_{xt}-2\,uu_{x^2}+2\,u_xu_{x^3}-u_{x^2}^2+\frac{4}{5}\,u_{x^3t}+2\,uu_{x^4},\\
\\
\textit{I}^{(t)}_{u^j}\textit{f}=u-\frac{1}{3}\,u_{x^2}+\frac{1}{5}\,u_{x^4},\end{array}\label{eq:8.14}\end{eqnarray}
apply (\ref{eq:8.9}) to the integrands (\ref{eq:8.14}), therefore
\begin{eqnarray} \begin{array}{lcl}
\Psi:={\textit{H}}^{(x)}_{u(x,t)}\textit{f}=\frac{3}{2}\,u^2-\frac{1}{2}\,u_{x}^2-\frac{2}{3}\,u_{xt}-uu_{x^2}+u_xu_{x^3}-\frac{1}{2}\,u_{x^2}^2
+\frac{4}{5}\,u_{x^3t}+uu_{x^4},\\
\\
\Phi:={\textit{H}}^{(t)}_{u(x,t)}\textit{f}=u-\frac{1}{3}\,u_{x^2}+\frac{1}{5}\,u_{x^4},
\end{array}\label{eq:8.15}\end{eqnarray}
so, we have the first conservation low of the higher order CH
equation respect to multiplier $\xi=1$
\begin{eqnarray} \begin{array}{lcl}
D_x\big(\frac{3}{2}\,u^2-\frac{1}{2}\,u_{x}^2-\frac{2}{3}\,u_{xt}-uu_{x^2}+u_xu_{x^3}-\frac{1}{2}\,u_{x^2}^2
+\frac{4}{5}\,u_{x^3t}+uu_{x^4}\big)+D_t\big(u-\frac{1}{3}\,u_{x^2}+\frac{1}{5}\,u_{x^4}\big)=0.
\end{array}\label{eq:8.16}\end{eqnarray}
Now we find conservation law respect to multiplier $\xi=u$, in
this case we have
\begin{eqnarray}
\textit{f}=U(U_t-U_{x^2t}+U_{x^4t}+3UU_x-2U_xU_{x^2}-UU_{x^3}+2U_{x}U_{x^4}+UU_{x^5}),\label{eq:8.17}\end{eqnarray}
the integrands are
\begin{eqnarray}
\begin{array}{lcl}
\textit{I}^{(x)}_{u^j}\textit{f}=3\,u^{3}+3\,u^2u_{x^4}-\frac{4}{3}\,uu_{xt}+\frac{2}{3}\,u_{x}u_{t}-3\,u^2u_{x^2}+\frac{8}{5}\,uu_{x^3t}
-\frac{2}{5}\,u_{x^3}u_{t}-\frac{6}{5}\,u_xu_{x^2t}+\frac{4}{5}\,u_{x^2}u_{xt},\\
\\
\textit{I}^{(t)}_{u^j}\textit{f}=u^2-\frac{2}{3}\,uu_{x^2}+\frac{1}{3}\,u_{x}^2+\frac{2}{5}\,uu_{x^4}-\frac{2}{5}\,u_{x}u_{x^3}
+\frac{1}{5}\,u_{x^2}^2,\end{array}\label{eq:8.18}\end{eqnarray}
applying 2-dimensional homotopy operator, we have
\begin{eqnarray} \begin{array}{lcl}
\Psi:={\textit{H}}^{(x)}_{u(x,t)}\textit{f}=u^3+u^2u_{x^4}-\frac{2}{3}\,uu_{xt}+\frac{1}{3}\,u_{x}u_{t}-u^2u_{x^2}+\frac{4}{5}uu_{x^3t}
-\frac{1}{5}\,u_{x^3}u_t-\frac{3}{5}\,u_xu_{x^2t}+\frac{2}{5}\,u_{x^2}u_{xt},\\
\\
\Phi:={\textit{H}}^{(t)}_{u(x,t)}\textit{f}=\frac{1}{2}\,u^2-\frac{1}{3}\,uu_{x^2}+\frac{1}{6}\,u_{x}^2+\frac{1}{5}\,uu_{x^4}
-\frac{1}{5}\,u_{x}u_{x^3}+\frac{1}{10}\,u_{x^2}^2,
\end{array}\label{eq:8.19}\end{eqnarray}
so, the second conservation low of the higher order CH equation is
\begin{eqnarray} \begin{array}{lcl}
D_x\big(u^3+u^2u_{x^4}-\frac{2}{3}\,uu_{xt}+\frac{1}{3}\,u_{x}u_{t}-u^2u_{x^2}+\frac{4}{5}uu_{x^3t}
-\frac{1}{5}\,u_{x^3}u_t-\frac{3}{5}\,u_xu_{x^2t}+\frac{2}{5}\,u_{x^2}u_{xt}\big)\\
\\
+D_t\big(\frac{1}{2}\,u^2-\frac{1}{3}\,uu_{x^2}+\frac{1}{6}\,u_{x}^2+\frac{1}{5}\,uu_{x^4}
-\frac{1}{5}\,u_{x}u_{x^3}+\frac{1}{10}\,u_{x^2}^2\big)=0.
\end{array}\label{eq:8.20}\end{eqnarray}

\section{Conclusion  }
In this paper by applying the criterion of invariance of the
equation under the infinitesimal prolonged infinitesimal
generators, we find the most general Lie point symmetries group
for higher order CH equation. By applying the nonclassical
symmetry method for the higher order CH equation, we concluded
that the analyzed model do not admit supplementary, nonclassical
type symmetries. Also, we have constructed the optimal system of
one-dimensional subalgebras for the higher order CH equation. The
latter, creates the preliminary classification of group invariant
solutions. The Lie invariants and similarity reduced equations
corresponding to infinitesimal symmetries are obtained. We find
the conservation laws from the multiplier method.

\end{document}